\documentclass[11pt]{article}

\usepackage{amssymb}\usepackage{amsfonts,amsmath,amscd}
\usepackage[T2A]{fontenc}
\usepackage[cp1251]{inputenc}
\usepackage[ukrainian,russian,english]{babel}

\begin{document}

\def\R{{\mathbb{R}}}
\def\C{{\mathbb{C}}}
\def\Z{{\mathbb{Z}}}
\def\N{{\mathbb{N}}}

\def\supp{{\rm supp\,}}
\def\intt{{\rm int\,}}
\def\a{\alpha}
\def\b{\beta}
\def\l{\lambda}
\def\d{\delta}
\def\e{\varepsilon}
\def\Im {\rm Im \;}
\def\Re {\rm Re \;}
\def\f {\varphi}

\large

\begin{center} Subharmonic Almost Periodic Functions
\end{center}

\begin{center}
Favorov S.Yu. Rakhnin A.
\end{center}

\begin{center}
V.N. Karazin Kharkiv National University,

4 Svobody sq., Kharkiv, 61077, Ukraine.

email: rakhnin@univer.kharkov.ua,
Sergey.Ju.Favorov@univer.kharkov.ua

\bigskip

\end{center}

\newpage

\bigskip

Subharmonic almost periodic functions

A.V. Rakhnin, S. Yu. Favorov

\textit{ We prove that almost periodicity in the sense of
distributions coincides with almost periodicity with respect to
Stepanov's metric for the class of subharmonic functions in a
horizontal strip. We also prove that Fourier coefficients of these
functions are continuous functions in  Im z. Further, if the
logarithm of a subharmonic almost periodic function is a
subharmonic function, then it is almost periodic.}

\bigskip

\newpage

Subharmonic almost periodic functions were introduced in
\cite{rrf} in connection with investigation of zero distribution
of holomorphic almost periodic functions in a strip. In this paper
almost periodicity was defined in the sense of distributions,
namely as almost periodicity of the convolution with a test
function. However, subharmonic functions $\log|f(z)|$, where
$f(z)$ is a holomorphic almost periodic function, were considered
much earlier in papers \cite{r1} and \cite{r2}, where the
important point was to prove almost periodicity of such functions
in the sense of distributions. In \cite{rakh} this was extended to
a subharmonic uniformly almost periodic function whose logarithm
is a subharmonic function.

In this paper we prove that subharmonic almost periodic in the
sense of distributions functions are almost periodic in the
classical sense, if we consider Stepanov integral metric instead
of the uniform metric. Therefore the classes of subharmonic almost
periodic in the sense of distributions functions and subharmonic
Stepanov almost periodic functions are the same.

Now the Fourier-Bohr coefficients of such functions can be defined
in the usual way. For a horizontal strip these coefficients are
functions depending on ${\rm Im} z$. In this paper we prove that
these coefficients depend continuously on ${\rm Im} z$, which
allows us to approximate any subharmonic almost periodic function
by exponential sums with continuous coefficients in Stepanov
metric. Thus we prove that subharmonic almost periodic functions
are Stepanov almost periodic in the sense of the definition in
\cite{udod}.

In \cite{rrf} it was proved that $\exp(u)$, where $u$ is a
subharmonic almost periodic in the sense of distributions
function, is also almost periodic in the sense of distributions.
Moreover, for an almost periodic function $\log|f(z)|$, where
$f(z)$ is a holomorphic function, $|f(z)|$ is uniformly almost
periodic. Conversely, we prove that the logarithm of a subharmonic
almost periodic function is an almost periodic function, provided
it is a subharmonic function. Thus we obtain a stronger than the
one in \cite{rakh}, as well as the converse to the result in
\cite{rrf}.

\bigskip

We start with the following definitions and notations (see
\cite[p. 51]{cord}).

{\bf Definition 1.} { \it A continuous function $f(z)$ ($z=x+iy$),
defined on ${\mathbb{R}}+iK$, where $K$ is a compact subset of
${\mathbb{R}}$ (it is allows that $K=\{0\}$), is called uniformly
almost periodic (Bohr almost periodic), if from any sequence ${\{
t_n\} \subset {\mathbb{R}}}$ one can choose a subsequence $\{
t_{n'}\}$ such that the functions $f(z+t_{n'})$ converge uniformly
on ${\mathbb{R}}+iK$.}

Equivalent definition is the following:

For any $\varepsilon>0$ there exists $L(\varepsilon)>0$ such that
each interval of length $L(\varepsilon)$ contains a real number
$\tau$ with the property
$$
\sup\limits_{z\in {\mathbb{R}} +i K}|f(z+\tau)-f(z)|<\varepsilon.
$$

{\bf Definition 2.} { \it A distribution $f(z)\in D'(S)$ of order
0 ($S$ is an open horizontal strip) is called almost periodic, if
for any test function $\f\in D(S)$ the convolution
$$
\int u(z)\f(z-t)dxdy
$$
is uniformly almost periodic on the real axis.}

Note that according to \cite{r2}, for an almost periodic
distribution $f(z)$ from any sequence $\{h_n\}\subset
{\mathbb{R}}$ one can choose a subsequence $\{h_{n'}\}$, such that
$\int \varphi(z)f(z+h_{n'})dxdy$ converge uniformly on every set
${\Gamma_{K} = \{ \varphi(z+t): t\in {\mathbb{R}}, \varphi\in K
\}}$, where $K$ -- is a compact subset of $D(S)$.

Any subharmonic function is locally integrable, so we can consider
it as a distribution.

A class of subharmonic almost periodic functions in an open strip
$S$ will be denoted by $WAP(S)$.

Furthermore, for $-\infty<\alpha<\b<+\infty$ we define
$$
S_{[\alpha,\beta]}=\{z\in{\mathbb{C}}: \alpha\le{\rm Im} z \le
\beta\},
$$
$$
{\rm Im} S = \{{\rm Im} z: z\in S\},
$$
and for functions $u$, $v$, which are integrable on horizontal
intervals in $S_{[\alpha,\beta]}$, we denote
$$
d_{[\alpha,\beta]}(u,v):=\sup_{z\in S_{[\alpha;\beta]}}\int_0^1
|u(z+t)-v(z+t)|dt.
$$

{\bf Definition 3.} { \it A function $f(z)$ integrable on
horizontal intervals in an open horizontal strip ${S}$ is called
Stepanov almost periodic, if from any sequence ${\{ h_n\} \subset
{\mathbb{R}}}$ one can choose a subsequence $\{ h_{n'}\}$ and a
function $g(z)$ such that the functions ${f(z+h_{n'})}$ converge
to $g(z)$ in the topology defined by seminorms
$d_{[\alpha,\beta]}$, $\alpha,\beta \in {\rm Im} S$}.

A class of a subharmonic Stepanov almost periodic functions in an
open strip $S$ will be denoted by $StAP(S)$. Since such functions
are Stepanov almost periodic on every line $y=const$, for $u\in
StAP(S)$ there exists the mean value
$$
M(u,y):=\lim_{T\to\infty}\frac{1}{2T} \int\limits_{-T}^{T}
u(x+iy)dx.
$$
To each such $u$ we can associate Fourier-Bohr series
$$
u(z)\sim\sum_{\l\in{\mathbb{R}}}a_{\l}(u,y)e^{i\lambda x},
$$
where
$$
a_{\l}(u,y):=M(u e^{-i\lambda x},y).
$$
are Fourier-Bohr coefficients.

{\bf Definition 4.} { \it A function $u(z)\ge 0$ is called
logarithmic subharmonic in a domain $G\subset {\mathbb{C}}$, if
the function $\log u(z)$ is subharmonic in this domain.}

It is easy to see that a logarithmic subharmonic function is
subharmonic.

\bigskip

We prove the following theorems:

\bigskip
{\bf Theorem 1.} { \it $u(z)\in WAP(S)$ if and only if $u(z)\in
StAP(S)$.}

\bigskip
{\bf Theorem 2.}  { \it Let $u(z)$ be a logarithmic subharmonic
function in a strip $S$. Then $\log u(z)\in WAP(S)$ if and only if
$u(z)\in WAP(S)$.}

\bigskip

{\bf Theorem 3.} {\it Let $u(z)$ be a subharmonic almost periodic
function in a strip $S$. Then its Fourier-Bohr coefficients are
continuous in ${\rm Im} S$.}

\bigskip

From Theorem 3 and Bessel inequality for Fourier-Bohr coefficients
it follows that spectrum of an almost periodic subharmonic
function $u(z)$ (i.e. the set $\{\l\in{\mathbb{R}}:
a_{\l}(u,y)\not\equiv 0\}$) it is most countable, which also
follows from Theorem~1.12 in \cite{r2}.

\bigskip

{\bf Theorem 4.} {\it Subharmonic function $u(z)$ in an open
horizontal strip $S$ is almost periodic if and only if there
exists a sequence finite exponential sums
\begin{eqnarray}
P_{m}(z)=\sum_{n=1}^{N_m} a_{n}^{(m)}(y) e^{i\lambda
_n^{(m)}x},\label{epol}
\end{eqnarray}
where $\lambda_n \in {\mathbb{R}}$, $a_{n}^{(m)}(y)\in C({\rm Im}
S)$, which converges to the function $u(z)$ in the topology
defined by seminorms $d_{[\alpha,\beta]}$, $\alpha,\beta \in {\rm
Im} S$.

Moreover, $P_{m}(z)$, $m=1,2,..$ are subharmonic functions in $S$.
}

\bigskip

To prove the theorems above we use the following propositions:

\bigskip
{\bf Proposition 1.} {\it Convergence of subharmonic functions in
$D'(G)$ is equivalent to the convergence in $L_{loc}^1(G)$}
(see~\cite{herm}).

\bigskip
{\bf Proposition 2.} {\it Weak limit of subharmonic functions is
subharmonic function} (see~\cite{herm}).

\bigskip

We denote by $G^{\mu}$ the Green potential of a measure $\mu$ for
the disk $B(R,0)$, i.e.
$$
G^{\mu}(z):=\int_{B(R,z_0)}\log \frac{|R^2 - z\overline \zeta|}
{R|z-\zeta|}d\mu(\zeta).
$$

\bigskip

{\bf Lemma 1.} { \it Let measures $\mu_n$ converge uniformly to a
measure $\mu$ in a neighborhood of the disk $\overline{B(R,0)}$,
and $\mu(\partial B(R,0))=0$. Then for any $t_1>0$, $t_2>0$ such
that $t_1^2+t_2^2<R^2$,
\begin{eqnarray}
\lim_{n\to\infty}\sup_{y\in
[-t_2;t_2]}\int\limits_{-t_1}^{t_1}\left|G^{\mu_n}(z)-G^{\mu}(z)\right|dx=0,
\label {sconv}
\end{eqnarray}
where $z=x+iy$. }

{P r o o f.}  Denote $\nu_n=\mu_n - \mu$. We have
$$
\sup_{y\in
[-t_2;t_2]}\int\limits_{-t_1}^{t_1}\left|G^{\mu_n}(z)-G^{\mu}(z)\right|dx\le
\sup_{y\in
[-t_2;t_2]}\int\limits_{-t_1}^{t_1}\left|\int_{B(R,0)}\log
\frac{|R^2 - z\overline \zeta |} {R} d\nu_n(\zeta)\right|dx+
$$
\begin{eqnarray}
+\sup_{y\in
[-t_2;t_2]}\int\limits_{-t_1}^{t_1}\left|\int_{B(R,0)}\log
|z-\zeta| d\nu_n(\zeta)\right|dx \label{4}
\end{eqnarray}

The condition $\mu(\partial B(R,0))=0$ implies that the
restrictions of the measures $\mu_n$ to the disk
$\overline{B(R,0)}$ converge weakly to the restriction of the
measure $\mu$ on the disk, and the function ${\log ({|R^2 -
z\overline \zeta|}R^{-1})}$ is continuous for $|x|\le t_1$,
$|y|\le t_2$, $\zeta \in \overline{B(R,0)}$. Thus the first term
on the right-hand side of (\ref{4}) is small. Without loss of
generality, we can assume that $R<1/2$, so that for ${z,\zeta \in
B(R,0)}$ we have $\log|z-\zeta|<0$.

Let $\varepsilon>0$ be an arbitrary fixed number. We denote
$\log_{\varepsilon} |z-\zeta|=\max\{\log|z-\zeta|,\log
\varepsilon\}$. This function is continuous for $|x|\le t_1$,
$|y|\le t_2$, $\zeta \in \overline{B(R,0)}$ and for any
$\varepsilon
>0$. We have

$$
\sup_{y\in
[-t_2;t_2]}\int\limits_{-t_1}^{t_1}\left|\int_{B(R,0)}\log
|z-\zeta| d\nu_n(\zeta)\right|dx\le \sup_{y\in
[-t_2;t_2]}\int\limits_{-t_1}^{t_1}\left|\int_{B(R,0)}\log_{\varepsilon}
|z-\zeta| d\nu_n(\zeta)\right|dx+
$$

$$
+\sup_{y\in [-t_2;t_2]}\int\limits_{-t_1}^{t_1}\int_{B(R,0)}\left
|\log |z-\zeta|-\log_{\varepsilon} |z-\zeta| \right |
d|\nu_n|(\zeta)dx.
$$
The first term on the right-hand side of this inequality is small
when $n$ is sufficiently large. Then
$$
\sup_{y\in [-t_2;t_2]}\int\limits_{-t_1}^{t_1}\int_{B(R,0)}\left
|\log |z-\zeta|-\log_{\varepsilon} |z-\zeta| \right |
d|\nu_n|(\zeta)dx=
$$
$$
\sup_{y\in [-t_2;t_2]}\int_{B(R,0)}\int\limits_{[-t_1;t_1]\cap
\{x: |x+iy-\zeta|\le \varepsilon\}}(\log \varepsilon - \log
|z-\zeta| ) dx d|\nu_n|(\zeta)\le
$$
$$
\int_{B(R,0)}\int\limits_{-\varepsilon}^{\varepsilon}(\log
\varepsilon-\log |x|) dx d|\nu_n|(\zeta)\le 2\varepsilon
|\nu_n|(B(R,0)).
$$
Note that since $\varepsilon>0$ is arbitrary, and measures $\nu_n$
weakly converge to zero, one can choose a constant
$C\in{\mathbb{R}}$ with ${|\nu_n|(B(R,0))<C}$. The lemma is
proved.

\bigskip

{\bf Lemma 2.} { \it Let $u_n(z)$ be a sequence of subharmonic
functions in a domain ${G \subset {\mathbb{C}}}$ converging to a
function $u_0(z)\not \equiv -\infty$ in $D'(G)$, and let
$$
\sup\limits_{z\in G'}u_n(z)\le W(G')<\infty
$$
for any subdomain $G'\subset G$. Then for any rectangle
${[a;b]\times [ \alpha ;\beta ]\subset G}$,}
\begin{eqnarray}
\lim_{n\to\infty}\sup_{y\in
[\alpha;\beta]}\int\limits_{a}^{b}\left|u_n(z)-u_0(z)\right|dx=0.
\label {nconv}
\end{eqnarray}

{P r o o f.} For every disk $B(z_0,R)\subset\subset G$ we have the
following representation
$$
u_n(z)=-G_{R}^{\mu_{n}}(z;z_0)+H_{R}(z;z_0;u_n) \;\; n=0,1.2..,
$$
where $\mu_n$ are the Riesz measures of the functions $u_n(z)$,
$G_{R}^{\mu_{n}}(z;z_0)$ is the Green potential of the measure
$\mu_n$ in the disk $B(z_0,R)$, and $H_R(z;z_0;u_n)$ are the best
harmonic majorants of the functions $u_n(z)$ in this disk.
Conditions of the lemma imply that $\mu_n$ converge weakly to the
measure $\mu_0$. Without loss of generality, we can assume that
$\mu(\partial B(z_0,R))=0$, and using Lemma~1 we conclude that for
any $t_1,t_2$, $t_1^2+t_2^2<R^2$

$$
\sup\limits_{y_0-t_2\le y \le
y_0+t_2}\int\limits_{-t_1+x_0}^{t_1+x_0}|G_{R}^{\mu_{n}}(z;z_0)-G_{R}^{\mu_{0}}(z;z_0)|dx\longrightarrow
0,
$$
when $n\to \infty$. From this it follows that the functions
$H_R(z;z_0;u_n)$ converge to the function $H_R(z;z_0;u_0)$ in
$D'(B(z_0,R))$. Now using the mean value property, Harnak
inequality, and obvious inequality
$$
H_R(z;z_0;u_n)\le W(B(z_0,R))<\infty, \;\; n=0,1,2..,
$$
we obtain uniform convergence of $H_R(z;z_0;u_n)$ to the function
$H_R(z;z_0;u_0)$ in the rectangle
$[-t_1+x_0,t_1+x_0]\times[-t_2+y_0,t_2+y_0]$. Covering the
rectangle $[\alpha,\beta]\times [a,b]$ by a finite number of such
rectangles, we prove the lemma.

{P r o o f of Theorem 1.} Inclusion $StAP(S) \subset WAP(S)$ is
obvious. We prove the opposite inclusion. We consider arbitrary
substrip $S_{[\alpha,\beta]}$, $\alpha,\beta \in {\rm Im} S$ and a
sequence $\{h_{j}\} \subset {\mathbb{R}}$. Since $u(z)$ is a
subharmonic almost periodic distribution, there exists a
subsequence $\{h_{j_k}\}$ such that for some subharmonic (clearly
also almost periodic) function $v(z)$ and for any $\varphi\in
D(S_{[\alpha,\beta]})$, uniformly in $t\in {\mathbb{R}}$,
\begin{eqnarray}
\lim \limits_{k\to \infty} \int_S (u(z+h_{j_k}+t) - v(z+t))\varphi
(z)dxdy=0.\label {usl}
\end{eqnarray}
Now we will show that the functions $u(z+h_{k})$ converge to
$v(z)$ in the topology defined by seminorms $d_{[\alpha,\beta]}$,
$\alpha,\beta \in {\rm Im} S$. Assuming the contrary, there exist
$\varepsilon_0>0$, $\alpha,\beta \in {\rm Im} S$ such that for an
infinite sequence $k'$
$$
d_{[\alpha,\beta]}(u(z+h_{j_{k'}}),v(z))>\varepsilon_0,
$$
and therefore there exists a subsequence $\{t_{k'}\} \in
{\mathbb{R}}$ such that
\begin{eqnarray}
\sup_{y\in [\alpha;\beta]}\int_0^1
|u(z+h_{j_{k'}}+t_{k'})-v(z+t_{k'})|dx>\varepsilon_0. \label {34}
\end{eqnarray}
Passing to a subsequence (if necessary), we can assume that
$$
u(z+h_{j_{k'}}+t_{k'})\to w(z) \; , \; v(z+t_{k'})\to w_1(z)\;
\hbox{в}\; D'(S_{[\alpha,\beta]}).
$$
Lemma~2 implies
$$
\sup\limits_{y\in
[\alpha,\beta]}\int\limits_{0}^{1}|u(z+h_{j_{k'}}+t_{k'})-w(z)|dx
\to 0, \;\; k'\to \infty
$$
and
$$
\sup\limits_{y\in
[\alpha,\beta]}\int\limits_{0}^{1}|v(z+t_{k'})-w_1(z)|dx \to 0,
\;\; k'\to \infty,
$$
and thus inequality (\ref{34}) implies
\begin{eqnarray}
\sup_{y\in [\alpha;\beta]}\int_0^1
|w(z)-w_1(z)|dx\ge\varepsilon_0.\label {contr1}
\end{eqnarray}
On the other hand, using (\ref {usl}), for any test function
$\varphi(z)$,
$$
\int_{S_{[\alpha,\beta]}} (w(z) - w_1(z))\varphi (z)dxdy= \lim
\limits_{k'\to \infty} \int_{S_{[\alpha,\beta]}}
(u(z+h_{j_{k'}}+t_k') - v(z+t_{k'}))\varphi (z)dxdy
$$
$$
=\lim \limits_{k'\to \infty} \int_{S_{[\alpha,\beta]}}
(u(z+h_{j_{k'}}) - v(z))\varphi (z-t_{k'})dxdy=0,
$$
and thus $w(z)=w_1(z)$ almost everywhere. Since $w(z)$ and
$w_1(z)$ are subharmonic functions, then $w(z) \equiv w_1(z)$,
which contradicts (\ref {contr1}). The theorem is proved.

\bigskip

To prove Theorem~2 we need the following lemmas.

{\bf Lemma 3.} {\it Let $\varphi (t)$ be a function continuous in
$[-c,c]$. Then for any $\varepsilon>0$ there exists $\delta$,
depending on $\varphi$ and $\varepsilon$, such that for two
integrable on compact set $K$ functions $f,g:\; K \to [-c,c]$ the
inequality
$$
\int\limits_K |f(x)-g(x)|dm<\delta
$$
implies the inequality
\begin{eqnarray}
\int\limits_K |\varphi (f(x))-\varphi
(g(x))|dm<\varepsilon\label{lem}.
\end{eqnarray}
}

{P r o o f.} Choose $\tau>0$ such that $|t_1-t_2|<\tau$ implies
$|\varphi(t_1)-\varphi(t_2)|<\frac{\varepsilon}{2m(K)}$, and
denote
$$
A_1=\{x\in K: |f(x)-g(x)|<\tau\},
$$
$$
A_2=\{x\in K: |f(x)-g(x)|\ge\tau\}.
$$
Notice that $m(A_2)\le \frac{1}{\tau} \int_{A_2} |f(x)-g(x)|dm$,
and therefore
$$
\int\limits_{K} |\varphi(f(x))-\varphi(g(x))|dm \le
\int\limits_{A_1}
|\varphi(f(x))-\varphi(g(x))|dm+\int\limits_{A_2}
|\varphi(f(x))-\varphi(g(x))|dm\le
$$
$$
\le \frac{m(A_1)\varepsilon}{2m(K)}+\frac{2 \sup |\varphi
(t)|}{\tau} \int\limits_{K}|f(x)-g(x)|dm.
$$
Choosing suitable $\delta$, (\ref{lem}) follows. The lemma is
proved.

\bigskip

{\bf Lemma 4.} { \it Let $u_n (z)$ be a sequence of uniformly
bounded from above logarithmic subharmonic functions in a domain
$G\subset {\mathbb{C}}$, converging to a function
$u_0(z)\not\equiv 0$ in the sense of distributions. Then the
functions $\log u_n (z)$ converge to the function $\log u_0(z)$ in
the sense of distributions.}

{P r o o f.} The functions $u_n(z)$ are logarithmic subharmonic,
and in particular subharmonic. Using Proposition~1, $u_n(z)$
converge to $u_0(z)$ in $L_{loc}^1(G)$.

Next, these functions are uniformly bounded from above by some
constant ${V>0}$, bounded from below by $0$, and the function
$l_{\varepsilon}(t)=\log \max\{\varepsilon,t\}$ is continuous in
the interval $[0,V]$. Lemma~3 implies that for fixed $\varepsilon$
the functions $l_{\varepsilon} (u_n)(z)$ converge to the function
$l_{\varepsilon} (u_0)(z)$ in $L_{loc}^1(G)$, and thus in the
sense of distributions. From Proposition~2 it follows that the
functions $l_{\varepsilon} (u_0)(z)$ are subharmonic for all
$\varepsilon$, and their monotone limit when $\varepsilon\to 0$,
i.e. the function $\log u_0$, is also subharmonic.

Now we consider a disk $B(z_0,r)\subset\subset G$. From the
convergence in $L^1 (B(z_0,r))$ of the sequence $u_n(z)$ it
follows that the subsequence $\{u_{n'}(z)\}$ converges uniformly
on every fixed compact set $K_1\subset B(z_0,r)$ with positive
Lebesgue measure. Since the function $\log u_0(z)$ is subharmonic
and not identically $-\infty$ on $K_1$,
$$
\sup\limits_{z\in K_1} (\log u_0(z))\ge C_0,
$$
or
$$
\sup\limits_{z\in K_1} (u_0(z))\ge e^{C_0}.
$$
Thus for all $n'>n_0$
$$
\sup\limits_{z\in K_1} (u_{n'}(z))\ge e^{C_0-1},
$$
and
$$
\sup\limits_{z\in K_1} (\log u_{n'}(z))\ge C_0-1, \; \forall
n=0,1..
$$
Since the functions $u_n(z)$ are uniformly bounded from above on
compact subsets of $G$, it follows that the family $\{\log u_{n'}
(z)\}$ is compact in $D'(G)$. Therefore there exists a subsequence
$\log u_{n''}(z)$ which converges in $D'(G)$ (and also in
$L_{loc}^1(G)$) to some subharmonic function $v(z)$ in $G$.

Note that for any compact set $K\subset G$ and for any
$\varepsilon>0$ we have the following inequality
$$
\int\limits_K |max\{\log u_{n'}(z),\log
\varepsilon\}-max\{v(z),\log \varepsilon\}|dxdy\le\int\limits_K
|\log u_{n'}(z)-v(z)|dxdy.
$$

Hence the functions $max\{\log u_{n'}(z),\log \varepsilon\}$
converge to the function $max\{v(z),\log \varepsilon\}$ in
$L^1_{loc}(G)$ for any $\varepsilon>0$.

On the other hand, as was shown above, $l_{\varepsilon}(u_n)(z)$
converge to $l_{\varepsilon} (u_0)(z)$ in $L^1_{loc}(G)$. Thus
almost everywhere (and, since the functions are subharmonic,
everywhere)

\begin{eqnarray}
\max \{v(z),\log\varepsilon\}=\max \{\log
u_0(z),\log\varepsilon\}. \label {1}
\end{eqnarray}

Since a set on which a subharmonic function equals to $-\infty$
has Lebesgue measure zero, then ${\varepsilon\to 0}$ implies that
${{\rm mes} (\{ z\in G: v(z)<\log \varepsilon\})\to 0}$, ${{\rm
mes} (\{ z\in G: \log u_0(z)<\log \varepsilon\})\to 0}$, and
${v(z)=\log u_0(z)}$ almost everywhere, and hence everywhere. Thus
the sequence of the functions $\log u_{n'}(z)$ converges to the
function $\log u_{0}(z)$ in $D'(G)$ and in $L^1_{loc}(G)$.

If for some subsequence of the functions $\log u_{n_j}(z)$,
$\varepsilon_0>0$ and compact set $K_0 \in G$
\begin{eqnarray} \int\limits_{K_0}
|\log u_{n_j}(z)-\log u_{0}(z)|dxdy\ge \varepsilon_0, \label
{prot}
\end{eqnarray}
then, using the above construction of the sequence $u_{n_j}(z)$,
we have that some subsequence of the sequence $\{\log u_{n_j}\}$
converges to $\log u_0 (z)$ in $L^1_{loc}(G)$, which contradicts
(\ref{prot}). The lemma is proved.

\bigskip

{P r o o f of Theorem 2.} From Proposition~3 in \cite{rrf} it
follows that the inclusion $\log u \in WAP(S)$ implies that
inclusion $u \in WAP(S)$. We are going to show the opposite
inclusion. Let $u(z)\in WAP(S)$ and ${\{h_n\}\subset
{\mathbb{R}}}$ be an arbitrary sequence. Passing to a subsequence
if necessary, we can assume that for some subharmonic function
$u_0$, uniformly in $t\in{\mathbb{R}}$,
\begin{eqnarray}
\lim \limits _{n\to \infty}\int_S (u(z+h_n)-u_0(z))\varphi (z-t)
dxdy=0.\label {uncon}
\end{eqnarray}
To prove the theorem it is sufficient to verify that uniformly in
$t\in{\mathbb{R}}$
\begin{eqnarray}
\lim_{n \to \infty} \int_S \log u(z+h_n+t)\varphi (z)dxdy=\int_S
\log u_0 (z+t)\varphi (z)dxdy. \label{*}
\end{eqnarray}
Assuming that this fails, for some $\varepsilon>0$ and some
sequence $t_n \to \infty$,
\begin{eqnarray}
\left | \int_S \log u(z+h_n+t_n)\varphi (z)dxdy-\int_S \log u_0
(z+t_n)\varphi (z)dxdy \right | \ge \varepsilon. \label {100}
\end{eqnarray}
Here $u_0(z)$ is a logarithmic subharmonic function with
${u_0(z)\in WAP(S)}$. Passing to a subsequence and using almost
periodicity of the function $u_0 (z)$, we can assume that
\begin{eqnarray}
\lim_{n \to \infty} \int_S u_0(z+t_n)\varphi (z)dxdy=\int_S v
(z)\varphi (z)dxdy \label{uncon2}
\end{eqnarray}
for some subharmonic in the strip S function $v(z)$. Since the
limit in~(\ref{uncon}) is uniform in $t\in {\mathbb{R}}$,
(\ref{uncon2}) implies
$$
\lim_{n \to \infty} \int_S u(z+h_n+t_n)\varphi (z)dxdy=\int_S v
(z)\varphi (z)dxdy.
$$
Now Lemma~4 implies that both integrals in~(\ref {100}) have the
same limit $\int \log v(z) \varphi (z) dxdy$, when $n \to \infty$,
which is impossible. Thus (\ref {*})~holds and Theorem~2 is
proved.

{P r o o f of Theorem 3.} Without loss of generality, we can
assume that $S$ is a strip with finite width. Let $S_0$ be an
arbitrary substrip, ${S_0}\subset\subset S$. Since the function
$u(z)$ is almost periodic, its Riesz measure $\mu:=\frac {1}{2\pi}
\Delta u $ is also almost periodic in the sense of distributions.

Denote
$$
K(w)=\frac{1}{2} \log | e^{-\gamma w^2} -1|,
$$
where
$$
0<\gamma<\frac {\pi}{\max\limits_{y_1,y_2\in {\rm Im}
S}(y_1-y_2)^2}.
$$

Note that the kernel $K(w)$ is a subharmonic function which is
bounded from above in $S$ and its restriction to $S_0$ satisfies
the equation
\begin{eqnarray}
\Delta K(w)=2\pi \delta (w), \label{delt}
\end{eqnarray}
where $\delta (w)$ is a standard Dirac measure. Denote
\begin{eqnarray}
V(z)=\int\limits_{S} K(w-z)\varphi({\rm Im} w)d\mu(w),\label{int}
\end{eqnarray}
where $\varphi\ge 0$ is a test function on ${\rm Im} S$ such that
$\varphi(y)=1$ for $y\in {\rm Im} S_0$.

Denote $P_n =\{(n-1/2,n+3/4) \times{\rm Im} S\} \subset S$. We are
going to show that $V(z)$ is a subharmonic function in every
$P_n$. Fixing $n_0\in{\mathbb{Z}}$, we have
$$
\int\limits_{S} K(w-z)\varphi({\rm Im}
w)d\mu(w)=\int\limits_{[n_0-1,n_0+1)\times{\rm Im} S}
K(w-z)\varphi({\rm Im} w)d\mu(w)+
$$

\begin{eqnarray}
+\sum_{n\in {\mathbb{Z}} \setminus
\{n_0-1,n_0\}}\int\limits_{[n,n+1)\times{\rm Im} S}
K(w-z)\varphi({\rm Im} w)d\mu(w).\label{row}
\end{eqnarray}
Every term in the right hand side of (\ref{row}) is obviously a
subharmonic function. For ${\rm Re} w \in [n,n+1)$, ${\rm Im} w
\in {\rm supp\,} \varphi$, $z\in P_{n_0}$, $n\ne n_0$, $n\ne
n_0-1$ we have
$$
\left | e^{-\gamma(w-z)^2} \right | = e^{-\gamma ({\rm Re} w -{\rm
Re} z)^2 + \gamma ({\rm Im} w - {\rm Im} z)^2}\le e^{\pi-\gamma
(|n-n_0|-3/4)^2}.
$$
Thus
$$
\sum_{n\in {\mathbb{Z}} \setminus \{n_0-1,n_0\}}\left
|\;\int\limits_{[n,n+1)\times{\rm Im} S} K(w-z)\varphi({\rm Im}
w)d\mu(w)\right | \le
$$
$$
\le \sum_{n\in {\mathbb{Z}} \setminus \{n_0-1,n_0\}}
\sup\limits_{z\in P_{n_0}} \sup \limits_{w\in [n,n+1)\times {\rm
supp\,} \varphi } \left |\frac {1}{2} \log |1-e^{-\gamma (z-w)^2}|
\right | \mu ([n,n+1)\times {\rm supp\,} \varphi).
$$
Since the measure $\mu$ is almost periodic, $\mu ([n,n+1)\times
{\rm supp\,} \varphi)$ is bounded from above uniformly in $n$ (see
\cite{rrf}), and therefore the series~(\ref{row}) converges
uniformly in $z\in P_{n_0}$ and the function $V(z)$ is subharmonic
in $P_{n_0}$, and also in $S$.

Now we are going to show that the function $V(z)$ is subharmonic
almost periodic in $S$. We consider a test function $\psi(z)$ on
$S$ and verify that the function
$$
f(t)=\int_{S} V(z) \psi (z-t) dz
$$
is uniformly almost periodic on the real axis. We have
$$
f(t)=\int_{S} \left (\int_{S} K(w-z)\psi (z) dz \right
)\varphi({\rm Im} w)d\mu(w+t).
$$
Note that the function
$$
\Psi (w):=\int_{S} K(w-z)\psi (z) dxdy
$$
is continuous in $S$, because the difference $K(w)-\log |w|$ is
continuous in some neighborhood of zero. Moreover, $\Psi
(w)=O(e^{-\gamma |w|^2})$ when $|{\rm Re} w| \to\infty$.

Since the values $\mu([n,n+1]\times {\rm Im} S_0)$ are uniformly
bounded in $n$, then
\begin{eqnarray} \label{mu}
\int \limits_{S} \frac {\varphi({\rm Im} w+{\rm Im}
z)d\mu(w+z)}{1+|w|^2}\le C_1<\infty
\end{eqnarray}
uniformly in $z\in S_0$.

We fix $\varepsilon>0$ and choose a test function $\nu(t)$, $0\le
\nu(t) \le 1$ on ${\mathbb{R}}$, and such that $\nu ({\rm Re}
w)=1$ on the set
$$
\left\{ w:\;|\Psi (w)|>\frac {\varepsilon}{C_1 (1+|w|^2)}\right\}.
$$
For all $t\in {\mathbb{R}}$ we have
$$
f(t)=\int_{S} \Psi (w) \nu({\rm Re} w) \varphi({\rm Im}
w)d\mu(w+t)+
$$
$$
+\int_{S} \Psi (w) (1-\nu({\rm Re} w)) \varphi({\rm Im}
w)d\mu(w+t).
$$
Property (\ref{mu}) implies that the second integral in the
equality is not greater than $\varepsilon$ for all $t\in
{\mathbb{R}}$. Since $\mu$ is an almost periodic measure, the
first integral is an almost periodic function, and if $\tau$ is an
$\varepsilon$-almost period, then it is a $2\varepsilon$-almost
period for $f$. Thus, the function $V(z)$ is a subharmonic almost
periodic, and in addition (\ref{delt}) implies that ${ \Delta V(z)
= 2\pi \varphi(y)\mu(z)}$ in the sense of distributions. Consider
the function
$$
H(z):=V(z) - u(z).
$$
This function is harmonic and almost periodic in the sense of
distributions in $S_0$. Let $\varphi\ge 0$ be a test function in
the disk $B(\varepsilon,0)$, which depends only on $|z|$ and such
that $\int \varphi (z) dxdy=1$. Since the convolution $\int H(z)
\varphi (z+\zeta) dxdy$ is equal to $H(\zeta)$ in some strip
$S_1\subset\subset S_0$, then the remark after Definition~2
implies uniform almost periodicity of the function $H(z)$ in
$S_1$. So its Fourier-Bohr coefficients are continuous in ${\rm
Im} S_1$ and, since $\varepsilon$ is arbitrary, in ${\rm Im} S_0$.
Thus it is enough show that the Fourier-Bohr coefficients of
$V(z)$,
$$
a_{\l}(V,y)=M(V e^{-i\lambda x},y)= \lim_{T\to\infty}\frac{1}{2T}
\int\limits_{-T}^{T} V(x+iy)e^{-i\lambda x}dx,
$$
are continuous.

We fix $\varepsilon>0$. We have
$$
K(w)=\max \{ K(w), -2\log N\} + \min \{ K(w) + 2\log N,
0\}=K_1(w)+K_2(w),
$$
where $N<\infty$ will be chosen later. Denote
$$
V_1(z):=\int\limits_{S} K_1(z-w)\varphi ({\rm Im} w) d\mu(w),
$$
$$
V_2(z):=\int\limits_{S} K_2(z-w)\varphi ({\rm Im} w) d\mu(w).
$$

Since for $|\gamma w^2|<1/2$ we have
$$
K(w)= 1/2 \log |1-1+\gamma w^2 -\frac {\gamma^2 w^4}{2!}
+...|=\log |w| + \theta (w),
$$
where $\theta (w)$ is a continuous function, then $K_2(w)=0$ for
$|w|\ge\delta_0>0$ and $N$ sufficiently large. Moreover, if
$|\theta (w)|\le \log N$, then for all $w\in {\mathbb{C}}$ and
$y\in {\rm Im} S$,
$$
\int\limits_{-T}^{T} K_2(z-w)dx = \int\limits_{-T}^{T}\min \{
\log|x-u+i(y-v)|+\theta (z-w) + 2\log N,0 \}dx \ge
$$
\begin{eqnarray}
\ge \int\limits_{-\infty}^{\infty}\min \{ \log|Nx-Nu|,0 \}dx
=-\frac{C}{N},\label{int1}
\end{eqnarray}
with some constant $C$, $0<C<\infty$. Now using the property that
${\mu ([n,n+1]\times {\rm supp\,} \varphi)}$ are bounded and the
fact that $K_2(z-w)=0$ for $|z-w|\ge \delta_0$, we have that for
all $T$
$$
\left | \frac{1}{2T} \int\limits_{-T}^{T} V_2(z) e^{-ix\l}
dx\right | \le \int\limits_{|{\rm Re} w|\le T+\delta_0}
\frac{1}{2T} \int\limits_{-T}^{T} |K_2(z-w)|dx d\mu(w)\le
$$
\begin{eqnarray}
\le\frac{C}{2TN} \int\limits_{|{\rm Re} w|\le T+\delta_0}
\varphi({\rm Im} w)d \mu(w)\le\frac{C_2}{N}\le
\varepsilon,\label{v2}
\end{eqnarray}
if $N$ is sufficiently large.

Further, since $K_1(w)=O(e^{-\gamma |w|^2})$ for $|{\rm Re} w|\to
\infty$, then one can choose a test function $0\le\eta(t)\le 1$ on
${\mathbb{R}}$ such that $\eta ({\rm Re} w)=1$ on the set
$$
\left \{ w:\;|K_1 (w)|> \frac {\varepsilon}
{C_1(1+|w|^2)}\right\},
$$
where $C_1$ is the constant from (\ref{mu}). We have
$$
V_1(z)=\int\limits_{S} K_1(w)\eta({\rm Re} w) \varphi ({\rm Im} w
+{\rm Im} z) d\mu(w+z)+
$$
$$
+\int\limits_{S} K_1(w) (1-\eta({\rm Re} w)) \varphi ({\rm Im} w
+{\rm Im} z) d\mu(w+z)=
$$
$$
=V_3(z)+V_4(z).
$$
From the choice of the function $\eta$ it follows that
\begin{eqnarray} \label{v4}
|V_4(z)|\le \varepsilon, \;\; \hbox{for} \;\; z\in S_0.
\end{eqnarray}
Since the kernel $K_1(w)$ is continuous and the family of shifts
of a test function in ${\rm Im} S_0$ is a compact set, then (see
the remark to Definition~2) the function $V_3(z)$ is uniformly
almost periodic in $S_0$ and it has continuous in ${\rm Im} S_0$
Fourier-Bohr coefficients (see \cite[p. 145]{cord}). Thus, if
${y_1,y_2 \in {\rm Im} S_0}$ and ${|y_1-y_2|\le \d(\varepsilon)}$,
then (\ref{v2}) and (\ref{v4}) imply
$$
|a_{\l}(V,y_1)-a_{\l}(V,y_2)|\le|a_{\l}(V_3,y_1)-a_{\l}(V_3,y_2)|+|a_{\l}(V_4,y_1)|
+|a_{\l}(V_4,y_2)| +
$$
$$
+|a_{\l}(V_2,y_1)|+|a_{\l}(V_2,y_1)|\le 5 \varepsilon .
$$
Thus $a_{\l}(V,y)$ are continuous. The theorem is proved.

{P r o o f of Theorem 4.} For $P_m(z)$ we choose Bohner-Fejer sums
of the function $u(z)$, i.e. set
$$
P_m (z):= \lim\limits_{T\to\infty} \frac
{1}{2T}\int\limits_{-T}^{T} u(z+t)\Phi^{(m)} (t) dt=\sum
k_{\l}^{(m)} a_{\l}(u,{\rm Im} z) e^{i\lambda {\rm Re} z}.
$$
Here $\Phi^{(m)} (t)$ is a sequence of Bohner-Fejer kernels (see.
\cite[p. 69]{lev}), and the set $\{ k_{\l}^{(m)}: k_{\l}^{(m)}\ne
0\}$ is finite for every $m$. Note that, according to Theorem~3,
the functions $a(u,y)$ are continuous in $y\in {\rm Im} S$.

We are going to show that $P_m(z)$ are subharmonic.

Note that the kernels $\Phi^{(m)} (t)$ are non-negative, bounded,
and
$$
\lim\limits_{T\to\infty} \frac {1}{2T}\int\limits_{-T}^{T}
\Phi^{(m)} (t) dt=1.
$$

Also note that the subharmonic almost periodic function $u(z)$ is
bounded from above in any subset $S'\subset\subset S$. Thus, using
Fatou's lemma, for any $m=1,2,...$, $z\in S$ and sufficiently
small $\rho$,
$$
\frac {1}{2\pi\rho}\int\limits_{0}^{2\pi} P_m (z+\rho
e^{i\varphi}) d\varphi \ge \overline{\lim\limits_{T\to\infty}}
\frac {1}{2T}\int \limits_{-T}^{T} \frac
{1}{2\pi\rho}\int\limits_{0}^{2\pi} u(z+\rho
e^{i\varphi}+t)\Phi^{(m)} (t)  d \varphi dt \ge
$$
$$
\ge \overline{\lim\limits_{T\to\infty}} \frac {1}{2T}\int
\limits_{-T}^{T} u(z+t)\Phi^{(m)} (t)dt=P_{m}(z).
$$

It is easy to see that the polynomials from (~\ref{epol}) are
almost periodic functions in the strip $S$ in the sense of
Definition~2, and therefore their limit in the topology defined by
seminorms $d_{[\alpha,\beta]}$, $\alpha,\beta \in {\rm Im} S$, is
an almost periodic function in the sense of Definition~2.

As it is shown in \cite{r2}, for any test function $\varphi(z)$ in
$S$ and for some (depending only on the spectrum of the function
$u(z)$) sequence of Bohner-Fejer sums, for $m\to \infty$,
uniformly in $t\in{\mathbb{R}}$

\begin{eqnarray} \label{sld}
\int\limits_{S} P_m(z)\varphi(z+t)dxdy \to \int\limits_{S}
u(z)\varphi(z+t)dxdy
\end{eqnarray}

Now we are going to verify that it implies the convergence in the
topology defined by seminorms $d_{[\alpha,\beta]}$, $\alpha,\beta
\in {\rm Im} S$. Indeed, if it is not true, then for some sequence
${x}_{m} \to \infty$ and some $\alpha,\beta\in {\rm Im} S$,
$\varepsilon_0>0$,
$$
\sup\limits_{y\in [\alpha,\beta]} \int_{0}^{1}
|u(x_m+iy+t)-P_m(x_m+iy+t)|dt\ge\varepsilon_0.
$$
Since $u\in StAP(S)$, then, passing to a subsequence if necessary,
one can assume that functions $u(z+x_m)$ converge to some function
$v\in StAP(S)$ with respect to metric $d_{[\alpha,\beta]}$, and
therefore
\begin{eqnarray} \label{svd2}
\sup\limits_{y\in [\alpha,\beta]} \int_{0}^{1}
|P_m(x_m+iy+t)-v(t+iy)|dt\ge\varepsilon_0/2.
\end{eqnarray}
Moreover, according to Theorem~2,
$$
\int\limits_{S} u(z+x_m)\varphi(z)dxdy \to \int\limits_{S}
v(z)\varphi(z)dxdy.
$$
Therefore, by setting $t=-x_m$ in (\ref{sld}), we have for any
test function $\varphi(z)$

$$
\lim \limits_{m\to\infty} \int\limits_{S} P_m(z+x_m)\varphi(z)dxdy
= \int\limits_{S} v(z)\varphi(z)dxdy.
$$
According to Lemma 2, this contradicts to (\ref{svd2}). The
theorem is proved.

\bigskip

\renewcommand{\refname}{\large \rm 

\end{document}